\input{style/arxiv-ba.cfg}
\documentclass[ba,linksfromyear,preprint]{imsart}
\makeatletter
   \@ifpackageloaded{natbib}{}{\usepackage{natbib}}
\makeatother
\usepackage{amssymb}
\usepackage{multirow}

\pubyear{2015}
\volume{10}
\issue{2}
\firstpage{275}
\lastpage{296}
\doi{10.1214/14-BA910}


\begin{document}

\begin{frontmatter}
\title{Dirichlet Process Hidden Markov Multiple Change-point Model}
\runtitle{Dirichlet Process Hidden Markov Multiple Change-point Model}

\begin{aug}
\author{\fnms{Stanley I. M.} \snm{Ko}\thanksref{t1}\ead[label=e1]{stanleyko@cuhk.edu.hk}},
\author{\fnms{Terence  T. L.} \snm{Chong}\thanksref{t2}\ead[label=e2]{chong2064@cuhk.edu.hk}},
\and
\author{\fnms{Pulak} \snm{Ghosh}\thanksref{t3}\ead[label=e3]{pulak.ghosh@iimb.ernet.in}}

\runauthor{S. I. M. Ko, T. T. L. Chong, and P. Ghosh}


\thankstext{t1}{Department of Finance and Business Economics, University of Macau, Macau,\hfill\break
{StanleyKo@umac.mo}}
\thankstext{t2}{Department of Economics and Institute of Global Economics and Finance, The Chinese University
of Hong Kong, Hong Kong, and Department of International Economics and Trade, Nanjing University, China,
{chong2064@cuhk.edu.hk}}
\thankstext{t3}{Department of Quantitative Methods \& Information Systems,
Indian Institute of Management at  Bangalore, India,
{pulak.ghosh@iimb.ernet.in}}

\end{aug}

\begin{abstract}
This paper proposes a new Bayesian multiple change-point
model which is based on the hidden Markov approach.
The Dirichlet process hidden Markov model does not require
the specification of the number of change-points \textit{a priori}.
Hence our model is robust to model specification
in contrast to the fully parametric Bayesian model.
We propose a general Markov chain Monte Carlo algorithm which only needs to
sample the states around change-points.
Simulations for a normal mean-shift model with known and unknown
variance demonstrate advantages of our approach. Two applications,
namely the coal-mining disaster data and the
real United States Gross Domestic Product growth, are provided. We detect a single
change-point for both the disaster data and US GDP growth.
All the change-point locations and posterior inferences of the two applications
are in line with existing methods.
\end{abstract}

\begin{keyword}
\kwd{Change-point}
\kwd{Dirichlet process}
\kwd{Hidden Markov model}
\kwd{Markov chain Monte Carlo}
\kwd{Nonparametric Bayesian}
\end{keyword}


\end{frontmatter}


\section{Introduction}\label{sec:intro}

The earliest Bayesian change-point model is explored by
\citet{chernoff_zacks64}, who assume a constant
probability of change at each point in time. \citet{smith75}
investigates the single change-point model under different
assumptions of model parameters. \citet{carlin_gelfand_smith92}
assume that the structural parameters are independent of the
change-points and introduce the Markov chain
Monte Carlo sampling method to derive the posterior distributions.
\citet{stephens94} further applies the Markov chain Monte Carlo
method to the case of multiple changes. \citet{chib98} allows the
change-point probability to depend on the regime between two adjacent
change-points. \citet{Koop_Potter2007} propose the Poisson hierarchical prior
for durations in the change-point model that allows the number of
change-points to be unknown. More recent works on the Bayesian change-point model
include \citet{wang_zivot2000}, \citet{giordani_kohn2006},
\citet{Pesaran_Pettenuzzo_Timmermann2006},
\citet{maheu_gordon08} and \citet{Geweke_yu2011}.

In this paper we follow the modeling strategy of \citet{chib98}
which is one of the most popular Bayesian change-point models.
He introduces a discrete
random variable indicating the regime from which a particular
observation is drawn. Specifically, let
$ Y_n = (y_1, y_2, \dots, y_n)' $ be the
observed time series, such that the density of $ y_t $ conditioned
on $ Y_{t-1} = ( y_1, y_2, \dots, y_{t-1} )' $ depends on the
parameter $ \theta $ whose value changes at an unknown time period
$ 1 < \tau_1 < \dots < \tau_k < n $ and remains constant within
each regime, that is,
\begin{equation}
  y_t \sim
  \begin{cases}
    \;p(y_t\mid Y_{t-1}, \theta_1) & \text{if $ t \leq \tau_1 $, }\\
    \;p(y_t\mid Y_{t-1}, \theta_2) & \text{if $ \tau_1 < t \leq \tau_2 $, }\\
    \;\vdots                         & \vdots                       \\
    \;p(y_t\mid Y_{t-1}, \theta_k)  & \text{if $ \tau_{k-1} < t \leq \tau_k $,}\\
    \;p(y_t\mid Y_{t-1}, \theta_{k+1}) & \text{if $ \tau_k < t \leq n $, }
  \end{cases}
\end{equation}
where $ \theta_i \in \mathbb{R}^l $ is an $l$ dimension vector, $ i = 1,
 2, \dots, k+1 $. Note that we consider in this paper the change-point
 problem when the data are assumed to be generated by a parametric model where
 the unknown parameter $ \theta_i $ changes  with respect to different regimes.
Let $ s_t $ be the discrete indicator variable such that
\begin{equation}
  y_t\mid s_t \sim p(y_t \mid Y_{t-1}, \theta_{s_t}),
\end{equation}
where $ s_t $ takes values in $ \{ 1, 2, \dots, k, k+1 \} $. The
indicator variable $ s_t $ is modeled as a discrete time,
discrete-state Markov process with the constrained transition
probability matrix
\begin{equation}\label{eq:tran_matrix}
  P =
  \begin{pmatrix}
    p_{11} & p_{12} & 0      & \cdots & 0\\
    0      & p_{22} & p_{23} & \cdots & 0\\
    \vdots & \ddots & \ddots & \ddots & \vdots\\
    \vdots & \vdots & \ddots & p_{kk} & p_{k(k+1)}\\
    0      & 0      & \cdots & 0      & 1
  \end{pmatrix},
\end{equation}
where $ p_{ij} = \mathrm{pr}( s_t = j \mid s_{t-1} = i ) $ is the
probability of moving to regime $ j $ at time $ t $ given that the
regime at time $ t - 1 $ is $ i $. With this parameterization, the
$ i $th change-point occurs at $ \tau_i $ if $ s_{\tau_i} = i $
and $ s_{\tau_i+1} = i + 1 $.

As pointed out in \citet{chib98}, the above is a hidden
Markov model where the transition matrix of the hidden state
$ s_t $ is restricted as in (\ref{eq:tran_matrix}).
Hence, Chib's multiple change-point model
inherits the limitation of the hidden Markov model in that the number
of states has to be specified in advance. In light of this,
\citet{chib98} suggests to select from alternative models (e.g. one
  change-point vs. multiple change-points) according to the Bayes
factors. In this paper, we introduce the Dirichlet process
hidden Markov model (DPHMM) with left-to-right transition dynamic, without
imposing restrictions on the number of hidden states. The use of the
DPHMM has the following appealing features:

\begin{enumerate}
\item We do not have to specify the number of states
\textit{a priori}. The information provided by the observations
determines the states endogenously. Hence, our method can be
regarded as semiparametric.
\item Our modeling approach facilitates the sampling of states since
we only need to sample the states around change-points.
\end{enumerate}
We note that \cite{kozumi2000} propose a method similar to ours
in that they utilize a Dirichlet process prior for
$\theta$, but in a mixture model.

The rest of the paper is organized as follows.
Section \ref{sec:dp} provides a brief introduction of the Dirichlet
process. Section \ref{sec:hmm_cpm} incorporates the Dirichlet
process into the change-point model. The general Markov chain
Monte Carlo sampler is discussed in Section~\ref{sec:MCMC}.
A~Monte Carlo study of the normal mean-shift model is conducted
in Section \ref{sec:normex}. Section \ref{sec:learn_albe} discusses learning
of DPHMM parameters. Section \ref{sec:app} provides
applications of our model and Section \ref{sec:rm}
concludes the paper.

\section{The Dirichlet Process}\label{sec:dp}
Our new method employs the Dirichlet process technique which is widely used in
nonparametric Bayesian models. The Dirichlet process prior is first proposed by
\citet{ferguson73}. He derives the Dirichlet process prior as the prior on the
unknown probability measure space with respect to some measurable space $
(\Omega, \mathcal{F}) $. Hence the Dirichlet process is a distribution over probability
measures. \citet{blackwell_macqueen1973} show that the Dirichlet process can be
represented by the Polya urn model. \citet{sethuraman1994} develops the
constructive sticking-breaking definition.

In the present study, we assume a Dirichlet process prior to each row of the
transition matrix. The Dirichlet process is best defined here as the
infinite limit of finite mixture models (\citet{neal1992}, \citet{neal2000}
and \citet{beal2002}). To illustrate the idea, let us first consider the case with
a finite number of states. With the left-to-right restriction to the
transition dynamic, a particular state $ s_{t-1} =i $ will either stay at the
current state $ i $ or transit to a state $ j > i $. A left-to-right Markov chain
with $ k $ states will typically have the following upper triangular
transition matrix
\begin{equation}\label{eq:ltr_tran_matrix}
  P =
  \begin{pmatrix}
    p_{11} & p_{12} & p_{13} & \cdots & p_{1k}\\
    0      & p_{22} & p_{23} & \cdots & p_{2k}\\
    0      & 0      & p_{33} & \cdots & p_{3k}\\
    \vdots & \vdots & \vdots & \ddots & \vdots\\
    0      & 0      & \cdots & \cdots & p_{kk}
  \end{pmatrix},
\end{equation}
where the summation of each row equals $ 1 $. Note that the
left-to-right Markov chain here is different from Chib's restricted band transition
matrix (\ref{eq:tran_matrix}). Here, the number of states $ k $ is not
necessarily the number of regimes as is the case in Chib's model.

Let $ \mathbf{p}_i = ( 0, \dots, p_{ii}, p_{i(i+1)}, \dots, p_{ik}) $
be the transition probabilities of the $ i $th row of the transition matrix
(\ref{eq:ltr_tran_matrix}). Suppose we draw $ m $ samples $ \{ c_1, \dots, c_m
\} $ of $ s_{t+1} $ given $ s_t = i $ with probability profile
$\mathbf{p}_i$. The joint distribution of the sample is thus
\begin{equation}
  \mathrm{pr}(c_1, \dots, c_m \mid  \mathbf{p}_i) = \prod_{j=i}^k p_{ij}^{m_j},
\end{equation}
where $ m_j $ denotes the number of samples that take
state $ j $, $ j = i, \dots, k $. We assume a symmetric Dirichlet prior $
\pi( \mathbf{p}_i\mid\beta) $ for $ \mathbf{p}_i $ with positive concentration
parameter $\beta $:\vadjust{\goodbreak}
\begin{equation}
  \mathbf{p}_i\mid\beta \sim \mathrm{Dirichlet}(\beta/(k-i+1),\dots,\beta/(k-i+1)) =
  \frac{\Gamma(\beta)}{\Gamma\left(\frac{\beta}{k-i+1}\right)^{k-i+1}}
  \prod_{j=i}^k p_{ij}^{\beta/(k-i+1)-1}.
\end{equation}
With the Dirichlet prior, we can analytically integrate out $ \mathbf{p}_i
$ such that
\begin{equation}
  \begin{aligned}
  \mathrm{pr}(c_1, \dots, c_m\mid\beta) &= \int \mathrm{pr}(c_1, \dots, c_m\mid
    \mathbf{p}_i,\beta) d{\bf \pi}(\mathbf{p}_i \mid\beta)\\
  &= \frac{\Gamma(\beta)}{\Gamma(m+\beta)} \prod_{j=i}^k \frac{\Gamma\left(m_j +
    \frac{\beta}{k-i+1}\right)}{\Gamma\left(\frac{\beta}{k-i+1}\right)}.
  \end{aligned}
\end{equation}
The conditional probability of a sample $ c_d \in \{ c_1, \dots,
c_m\} $ given all other samples is thus
\begin{equation}\label{eq:finit_state_cpmf}
  \mathrm{pr}(c_d = j \mid {\bf c}_{-d}) =
    \frac{m_{-d,j} + \beta/(k+i-1)}{m-1+\beta},
\end{equation}
where $ {\bf c}_{-d} $ denotes the sample set with $ c_d $
deleted, and $ m_{-d,j} $ is the number of samples in $ {\bf c}_{-d} $  that take
state $ j $.

Taking the limit of equation (\ref{eq:finit_state_cpmf}) as $ k $ tends
to infinity, we have
\begin{equation}\label{eq:finite_state_dpprior}
  \mathrm{pr}(c_d = j\mid{\bf c}_{-d}) =
  \begin{cases}
    \frac{m_{-d,j}}{m-1+\beta} & j\in \{ i, i+1, \dots, k \},\\
    \frac{\beta}{m-1+\beta}    & \text{for all potential states }.\\
  \end{cases}
\end{equation}
Note that the probability that $ c_d $ takes an existing state, say $
j $,  is proportional to $ m_{-d,j} $, which implies that $ c_d $ is more likely to
choose an already popular state. In addition, the probability
that a new state (i.e. $ k+1 $) takes place is proportional to $ \beta $. Hence, there
are potentially many states available, with infinite dimension transition
matrix
\begin{equation}\label{eq:ltr_tran_matrix2}
  P =
  \begin{pmatrix}
    p_{11} & p_{12} & p_{13} & \cdots \\
    0      & p_{22} & p_{23} & \cdots \\
    0      & 0      & p_{33} & \cdots \\
    \vdots & \vdots & \vdots & \ddots \\
  \end{pmatrix}.
\end{equation}
The actual state space can be regarded
as consisting of an infinite number of states, only a finite number of which
are actually associated with the data. Therefore, the number of states
is endogenously determined.

\section{The Dirichlet Process Hidden Markov Multiple Change-point Model and the State Evolution}\label{sec:hmm_cpm}

Let us now turn to the proposed multiple change-point model and discuss a particular
state evolution. Suppose we have already
generated the hidden states up to $ s_t = i $. We impose the Dirichlet process
as described in Section \ref{sec:dp} to $ s_{t+1} $. In the
change-point model, the transitions that have existed so far from state $ i $
are only self transitions. With the left-to-right restriction, we will
neither see a backward transition, i.e., transition from state $ i $ to some
previous states, nor a forward transition, i.e., transition to some future
states. The counts of the existing transitions from state $ i $ will be used as
the counts defined in equation (\ref{eq:finite_state_dpprior}). Hence, we will have
\begin{equation}\label{eq:prelim_state_prior}
  \mathrm{pr}(s_{t+1} = j\mid s_t = i, s_1, \dots, s_{t-1}) =
  \begin{cases}
    \frac{n_{ii}}{n_{ii}+\beta} & j = i,\\
    \frac{\beta}{n_{ii}+\beta}  & \text{ $ s_{t+1} $ takes a new state},\\
  \end{cases}
\end{equation}
where $ n_{ii} = \sum_{t'=1}^{t-1}\delta(s_{t'}, i)\delta(s_{t'+1},
i) $ denotes the counts of transitions that have occurred so far from state $ i
$ to itself.\footnote{The Kronecker-delta function $ \delta(a,b) = 1 $ if and
  only if $ a = b $ and $ 0 $ otherwise.}
Note that in equation (\ref{eq:prelim_state_prior}), $ s_{t+1} $ depends
only on the \textit{state} that $ s_t $ takes according to the
Markovian property. All other previous states merely provide the transition counts.

We introduce a self-transition prior mass $ \alpha
$ for each state. The idea here is that if $ s_{t} $  transits to a new
state, say $ s_{t+1} = i+1 $, then without previous transition records, the next state $
s_{t+2} $ conditioned on $ s_{t+1} = i+1 $ will further take another new
state with probability $ 1 $. Hence, with $ \alpha $, the trivial case is avoided and we have
\begin{equation}\label{eq:state_prior}
    \mathrm{pr}(s_{t+1} = j\mid s_t = i, s_1, \dots, s_{t-1}) =
  \begin{cases}
    \frac{n_{ii} + \alpha}{n_{ii}+\beta + \alpha} & j = i,\\
    \frac{\beta}{n_{ii}+\beta + \alpha}  & \text{ $ s_{t+1} $ takes a new state}.\\
  \end{cases}
\end{equation}

\begin{figure}[t]
\includegraphics{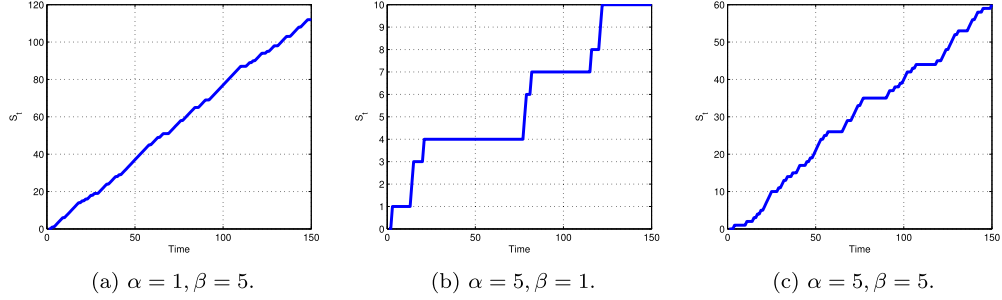}
\caption{Left-to-right Markov chain with different $\alpha$ and $\beta$.}
\label{fig:LtoR_mc}
\end{figure}

Therefore, the whole Markov chain is characterized by two parameters, $ \alpha
$ and $ \beta $, instead of a transition probability matrix. We can see that $
\alpha $ controls the prior tendency to linger in a state, and $ \beta $
controls the tendency to explore new states. Figure~\ref{fig:LtoR_mc} illustrates
three left-to-right Markov chains of length $ n = 150 $ with different $
\alpha $ and $ \beta $. Figure~\ref{fig:LtoR_mc}a depicts the chain that
explores many new states with very short linger time. Figure~\ref{fig:LtoR_mc}b
shows the chain with long linger time and less states.
Figure~\ref{fig:LtoR_mc}c lies in between.

Equation (\ref{eq:state_prior}) coincides with \citet{chib98}'s model when the
probability $p_{ii}$ is integrated out. Specifically, in \citet{chib98},
\begin{equation}\label{eq:Chib_int_p}
\begin{aligned}
        \mathrm{pr}(s_{t+1} = j&\mid s_t = i, s_1, \dots, s_{t-1}) \\
        &= \int p(s_{t+1} = j\mid s_t = i, s_1, \dots, s_{t-1}, p_{ii})
        f(p_{ii}) d p_{ii}\\ &=
  \begin{cases}
    \frac{n_{ii} + a}{n_{ii}+a + b} & j = i,\\
    \frac{b}{n_{ii}+a + b}  & j = i+1\\
  \end{cases}
\end{aligned}
\end{equation}
where $p_{ii} \sim \mathrm{Beta}(a,b)$ and $f(p_{ii})$ is the corresponding
density. However, our model stems from a
different perspective. The derivations in the previous section and equation
(\ref{eq:state_prior}) follow the nonparametric Bayesian literature
(\citet{neal2000}) and the infinite HMM of \citet{beal2002}. Indeed, it is
known that when the Dirichlet process is truncated at a finite number of
states, the process reduces to the generalized Dirichlet distribution
(GDD), see \citet{connor_mosimann1969}  and \citet{wong_GDD_1998}. For
the same reason we have (\ref{eq:state_prior}) coincides with
(\ref{eq:Chib_int_p}). We would like to point out that
our modeling strategy facilitates the Gibbs sampler of $ s_t $ which is
different from \citet{chib98}. We will elaborate further on the Gibbs sampler in the
next section. Also, learning of $\alpha$ and $\beta$ will be discussed in
Section \ref{sec:learn_albe}.\vspace*{-9pt}

\section{Markov Chain Monte Carlo Algorithm}\label{sec:MCMC}

\subsection{General}
Suppose we have observations $ Y_n =
(y_1, y_2, \dots, y_n )' $. Given the state $ S_n = (s_1, \dots, s_n)'
$ we have
\begin{equation}
  y_t\mid s_t \sim p(y_t \mid Y_{t-1}, \theta_{s_t}),
\end{equation}
where $ \theta_{s_t} \in \mathbb{R}^l $, $ Y_{t-1} =
 ( y_1, \dots, y_{t-1} )' $. Let $ \theta = (\theta_1,\dots,
\theta_k)' $ and $\gamma$ denotes a hyperparameter. Recall that we impose the
 DPHMM to the states
and a hierarchical model to the parameter,
we are thus interested in sampling from the
posterior $ p(\theta, S_n, \gamma \mid Y_n) $ given the priors
$ p(\theta|\gamma) $, $ p(\gamma) $ and $ p(S_n) $. The general Gibbs
sampler procedure is to sample the following in turn:

\vspace{3 mm}
\emph{Step} 1. $ S_n \mid \theta, \gamma, Y_n $,

\vspace{3 mm}
\emph{Step} 2. $ \theta \mid \gamma, S_n, Y_n $,

\vspace{3 mm}
\emph{Step} 3. $ \gamma \mid \theta, S_n, Y_n $.
\vspace{3 mm}

\noindent We will discuss the three steps below.

\subsection{Simulation of $ S_n $}

The state prior $ p(S_n) $ can be easily derived from
(\ref{eq:state_prior}). Moreover, the full conditional is
\begin{equation}\label{eq:post_state}
  p(S_n \mid \theta,\gamma, Y_n) \propto p(S_n) p(
    Y_n\mid S_n, \theta, \gamma).
\end{equation}
Simulation of $ S_n $ from the full conditional (\ref{eq:post_state}) is done
by the Gibbs sampler. Specifically, we draw $ s_t $ in turn for $ t = 1, 2,
\dots, n $ from
\begin{equation}\label{eq:gibbs_post_state}
  p(s_t \mid S_{t-1}, S^{t+1}, \theta, \gamma, Y_n) \propto
  p(s_t \mid s_{t-1}, S_{t-2}) p(s_{t+1} \mid s_t, S^{t+2}) p(y_t \mid s_t,
  Y_{t-1},\theta,\gamma),
\end{equation}
where $ S_{t-1} = ( s_1, \dots, s_{t-1} )' $ and
$ S^{t+1} = ( s_{t+1}, \dots, s_n )' $. The most recent
updated values of the conditioning variables are used in each iteration. Note that we
write $ p(s_t \mid s_{t-1}, S_{t-2}) $ and $  p(s_{t+1} \mid s_t,
  S^{t+2}) $ to emphasize the Markov dynamic; the other conditioning states
merely provide the counts.

With the left-to-right characteristic of the chain, we do not have to
sample all $s_t$ from $t=1$ to $T$. Instead, we only need to sample the state
in which a change-point takes place. To see this,
let us consider a concrete example. Suppose from the last
sampler, we have

\begin{center}
  \begin{tabular}{cccccccccccc}
    \hline
    $ s_1 $ & $ s_2 $ & $ s_3 $ &  $ s_4 $ &  $ s_5 $ &  $ s_6 $ & $ s_7 $ &
    $ s_8 $ & $ s_9 $ & $ s_{10}$ &  $ s_{11} $ & $ \cdots $ \\
    \hline
       1    &    1    &    1    &    2     &    2     &     3    &    3    &
       3    &    4    &    5    &    5     & $ \cdots $ \\
    \hline
  \end{tabular}
\end{center}

\noindent With left-to-right transition restriction, $ s_t $ requires a
sampling from (\ref{eq:gibbs_post_state}) if and only if $ s_{t-1} $ and $ s_{t+1} $
differ. For other cases, $ s_t $ is unchanged with probability one. Suppose we
are at $ t = 2 $. Since $ s_1 $ and $ s_3 $ are both equal to $ 1 $, $ s_2 $ is
forced to take $ 1 $. In the above chain,
the first state that needs to sample from (\ref{eq:gibbs_post_state}) is $ s_3
$, which will either take the values $ 1 $ or $ 2 $ ($s_2$ or $s_4$). If $ s_3 $
takes 2 in the sampling (i.e., joining the following regime), then the next state to sample
would be $ s_5 $; otherwise (i.e., joining the preceding regime), the next state
to sample is $ s_4 $ because $ s_5 - s_3 \neq 0 $ for $ s_3 = 1 $.
Now suppose we are at $ t = 9 $. We will draw a new $ s_9 $ and $ s_9 $ will
either join the regime of $ s_8 $ or the regime of $ s_{10} $. This will look
strange because a gap exists in the chain. However, our concern here is
the \emph{consecutive grouping} or \emph{clustering} in the series. We can alternatively
think that the state represented by $s_9$ is simply pruned away
in the current sweep. Note the numbers
assigned to the $ s_t $'s are nothing but indicators of regimes.\footnote{This
will exclude the case of the regime with only one data point. We do not
consider this situation here.} Therefore, we will relabel the $ s_t $'s after
each sweep.

In general, suppose $ s_{t-1} = i $ and $ s_{t+1} = i+1 $.
$ s_t $ takes either $ i $ or $ i+1 $. Table~\ref{tab:sampling_prob_state} shows
the corresponding probability values specified in (\ref{eq:gibbs_post_state}). To
see this, if $ s_t $ takes $ i $, then the transition from $ s_{t-1} $ to $
s_t $ is a self-transition and that from $ s_t $ to $ s_{t+1} $ is an
innovation. The corresponding probability values are in the first row of
Table~\ref{tab:sampling_prob_state}. The reasoning for
$ s_t = i+1 $ is similar. Note the changes of counts in
different situations.\vspace*{-6pt}

\begin{table}[htbp]\centering
\caption{Sampling probabilities of $ s_t $. $ n_{ii}
= \sum_{t'=1}^{t-2}\delta(s_{t'},i)\delta(s_{t'+1},i) $ and $ n_{i+1,i+1}
= \sum_{t'=t+1}^{n-1}\delta(s_{t'},i+1)\delta(s_{t'+1},i+1) $.}
\label{tab:sampling_prob_state}
\begin{tabular*}{\columnwidth}{@{\extracolsep{\fill}}lccc}
  \hline
  & $ p(s_t \mid s_{t-1} = i, S_{t-2})  $  &
  $ p(s_{t+1} = i+1 \mid s_t, S^{t+2}) $ &
  $ p(y_t\mid s_t,Y_{t-1},\theta) $  \\
  \hline
  \multirow{2}{*}{$ s_t = i $} &
  \multirow{2}{*}{$ \displaystyle \frac{n_{ii} + \alpha}{n_{ii}+\beta + \alpha} $} &
  \multirow{2}{*}{$ \displaystyle \frac{\beta}{n_{ii}+1+\beta + \alpha} $}  &
  \multirow{2}{*}{$ \displaystyle p(y_t\mid Y_{t-1},\theta_i) $} \\
  &&&\\
  \multirow{2}{*}{$ s_t = i+1 $}  &
  \multirow{2}{*}{$ \displaystyle \frac{\beta}{n_{ii}+\beta + \alpha} $}  &
  \multirow{2}{*}{$ \displaystyle \frac{n_{i+1,i+1} + \alpha}{n_{i+1,i+1}+\beta + \alpha} $}  &
  \multirow{2}{*}{$ \displaystyle p(y_t\mid Y_{t-1}, \theta_{i+1}) $} \\
  &&&\\
  \hline
\end{tabular*}
\end{table}
\noindent For the initial point $ s_1 $, if currently $ s_2 - s_1 \neq 0$,
then we can sample $ s_1 $ from
\begin{equation}
  \mathrm{pr}(s_1 \mid s_2, s_3, \dots, s_n) =
  \begin{cases}
    c\cdot\frac{\alpha}{\beta + \alpha}\cdot\frac{\beta}{\beta+\alpha}\cdot
    p(y_1\mid Y_0,\theta_{s_1})& \text{if $ s_1 $ unchanged, }\\
    c\cdot\frac{\beta}{\beta + \alpha}\cdot\frac{n_{s_2 s_2}+\alpha}{n_{s_2
        s_2}+\beta+\alpha}\cdot p(y_1\mid Y_0,\theta_{s_2})
    & \text{if $ s_1 = s_2 $, }\\
  \end{cases}
\end{equation}
where $ n_{s_2 s_2} = \sum_{t'=2}^{n-1}\delta(s_{t'},s_2)\delta(s_{t'+1},s_2)
$ and $ c $ is the normalizing constant. For the end point $ s_n $, if $
s_n- s_{n-1} \neq 1 $, we sample $ s_n $ from
\begin{equation}
  \mathrm{pr}(s_n \mid s_{n-1}, s_{n-2}, \dots, s_1) =
  \begin{cases}
    c\cdot\frac{n_{s_{n-1} s_{n-1}} + \alpha}{n_{s_{n-1} s_{n-1}} + \beta + \alpha}
        \cdot p(y_n\mid Y_{n-1},\theta_{s_{n-1}})& \text{if $ s_n = s_{n-1}$, }\\
    c\cdot\frac{\beta}{n_{s_{n-1} s_{n-1}} + \beta + \alpha}
      \cdot p(y_n\mid Y_{n-1},\theta_{s_n})& \text{if $ s_n $ unchanged, }\\
  \end{cases}
\end{equation}
where $ n_{s_{n-1} s_{n-1}} =
\sum_{t'=1}^{n-2}\delta(s_{t'},s_{n-1})\delta(s_{t'+1},s_{n-1}) $ and $ c $ is
the normalizing constant.

As mentioned in Section \ref{sec:hmm_cpm}, the DPHMM facilitates the Gibbs
sampler of states. In sampling $s_t$ from (\ref{eq:gibbs_post_state}),
we simultaneously use up all information of the transitions prior to $t$
(i.e. $s_1,\dots, s_{t-1}$) and after $t$ (i.e. $s_{t+1},\dots, s_T$) which
are captured in  $p(s_t \mid s_{t-1}, S_{t-2})$ and $p(s_{t+1} \mid s_t, S^{t+2})$.
Thus far, our algorithm only requires a record of transitions and draws
at the point where the structural change takes place, whereas we have to sample all
$s_t$ in \citet{chib98}.

\subsection{Updating $ \theta $ and $\gamma$}

Given $ S_n $ and $ Y_n $, the full conditionals of $ \theta $ and $\gamma$
are simply
\begin{equation}
  \begin{aligned}
    p(\theta_i \mid \gamma, S_n, Y_n) &\propto p(\theta_i \mid\gamma)\prod_{\{t:s_t=i\}}
    p(y_t \mid Y_{t-1}, \theta_i),\\
    p(\gamma \mid \theta, S_n, Y_n) &\propto p(\gamma) \;p(\theta \mid \gamma, S_n, Y_n),
  \end{aligned}
\end{equation}
which are model specific. In the following sections, we will
study a simulated normal mean-shift model, a discrete type Poisson
model, and an \textsc{ar}$(2)$ model.

\subsection{Initialization of States}
In Section 4.2, we have discussed the simulation of the states $S_n$. The
number of change-points is inherently estimated through the sampling of states
in equation (\ref{eq:gibbs_post_state}). Within the burn-in period, the state number
will be changing around after each MCMC pass.
After the burn-in period, the Markov chain converges and hence the number of states
becomes stable. Theoretically, it is legitimate to set any number of
change-points in the beginning and let the algorithm find out the convergent
number of states. In practice, we find it is more efficient to  initialize
with a large number of states and let the algorithm prune away redundant states,
rather than allow for the change-point number to grow from a small
number. Specifically, suppose a reasonably large state number $ k $ is
proposed. We initialize equidistant states, that is
\begin{equation}
  s_t = i, \quad  \mbox{if}\;\;  \frac{(i-1) \cdot n}{k} < t \leq  \frac{i \cdot n}{k},
\end{equation}
where $ i = 1,\dots, k$. Then the algorithm described above will work out the change-point locations
and the number of states after convergence of the Markov chain.

\section{A Monte Carlo Study: the Normal Mean-Shift Model}\label{sec:normex}

\subsection{The Model}

In this section, we first study the normal mean-shift model with known variance $
\sigma^2 $. Suppose the normal data $ Y_n = (y_1, \dots, y_n)' $ is
subject to unknown $ k $ changes in mean. We use the following
hierarchical model
\begin{equation}\label{eq:norm_model}
  \begin{aligned}
    &y_t\mid \theta_i \sim \mathrm{N}(\theta_i,
    \sigma^2)\;\;\;\;\text{if $ \tau_{i-1} < t \leq \tau_{i} $},\\
    &\theta_i\mid \mu, \upsilon^2 \sim \mathrm{N}(\mu, \upsilon^2),\\
    &(\mu, \upsilon^2) \sim \mbox{Inv-Gamma}(\upsilon^2\mid a, b),
  \end{aligned}
\end{equation}
where $ \sigma^2$ is known and $ \tau_i $ ($ i = 1, \dots, k $) is
the change-point. We set $ \tau_0 = 0 $ and $
\tau_{k+1} = n $. Next, we apply our algorithm to the case of unknown
variance. In addition to (\ref{eq:norm_model}), we assume the Inverse-Gamma
prior for $ \sigma^2 $
\begin{equation}\label{eq:sig2_prior}
  \sigma^2 \sim \mbox{Inv-Gamma}(\sigma^2\mid c,d).
\end{equation}
All derivations of the full conditionals and the Gibbs samplers are given in the Appendix.

\begin{figure}[b!]
\includegraphics{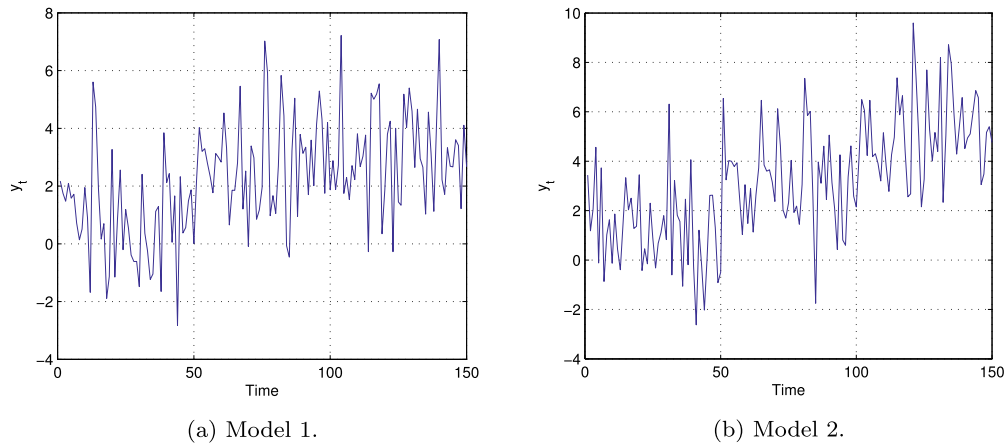}
\caption{Random realizations of Model 1 and Model 2.}
\label{fig:norm_realiz}
\end{figure}

We simulate two normal sequences with the parameters
specified in Table~\ref{tab:normodel1_2}. Specifically, Model 1 is subject to
one change-point occurring at $ t = 50 $. Model 2 is a two change-points model
with breaks at $ t = 50 $ and $ t = 100 $. Both models assume
variance $ \sigma^2 = 3 $. Two realizations with respect to Model 1 and Model 2
are shown in Figure~\ref{fig:norm_realiz}. We can see the overlapping of the
data ranges of different regimes and it is hard to visually identify the
change-points.

\begin{table}[htbp]
\caption{Normal Mean-Shift Models 1 and 2.}
\label{tab:normodel1_2}
\begin{tabular*}{\columnwidth}{@{\extracolsep{\fill}}lcccccccccccccccc}
  \hline
  && $ \theta_1 $ && $ \theta_2 $ && $ \theta_3 $ && $ \sigma^2 $
  && $ \tau_1 $   && $ \tau_2 $   && $ k $        && $ n $ \\
  \hline
      {\bf Model 1}      && \multirow{2}{*}{1} && \multirow{2}{*}{3} &&
      \multirow{2}{*}{-} && \multirow{2}{*}{3} && \multirow{2}{*}{50} &&
      \multirow{2}{*}{-} && \multirow{2}{*}{1} && \multirow{2}{*}{150}\\
      (One change-point) && &&   &&   &&   &&   &&    &&     &&    \\
      {\bf Model 2}      && \multirow{2}{*}{1} && \multirow{2}{*}{3} &&
      \multirow{2}{*}{5} && \multirow{2}{*}{3} && \multirow{2}{*}{50} &&
      \multirow{2}{*}{100} && \multirow{2}{*}{2} && \multirow{2}{*}{150}\\
      (Two change-points)&& &&   &&   &&   &&   &&    &&     &&       \\
      \hline
\end{tabular*}
\end{table}

\subsection{Simulation Results}

\begin{figure}
\includegraphics[scale=0.97]{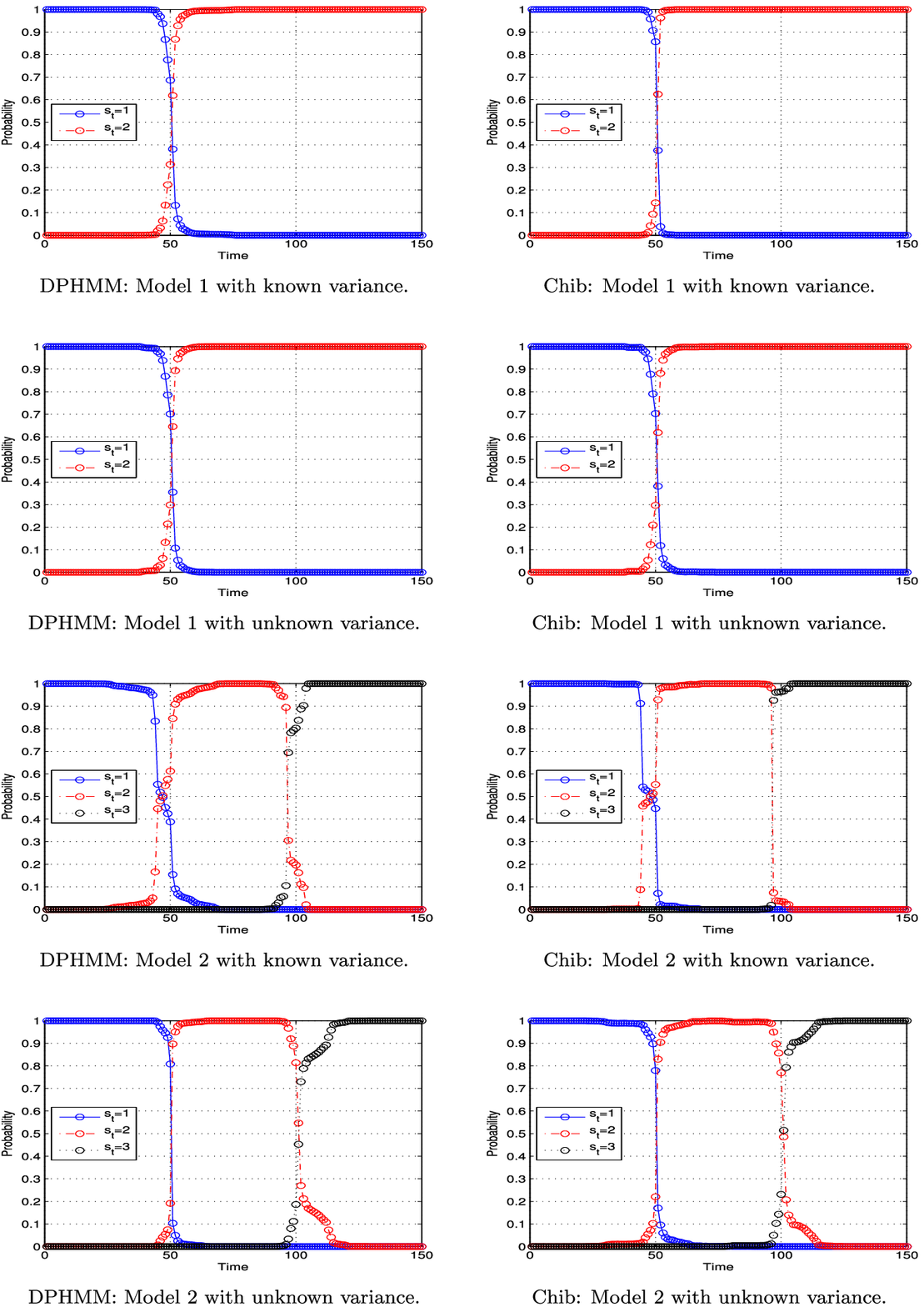}
\caption{Posterior probability of $ s_t = i $: DPHMM vs. Chib's model.}
\label{fig:normexhier_s}
\end{figure}

To implement our algorithm, we set the inverse-Gamma hyperparameters $ a = b =
c = d = 1 $, and the DPHMM parameters $ \alpha = 3 $ and $ \beta = 2 $.
The two Gibbs samplers for the cases of known and unknown variance are
conducted for 5000 sweeps with 5000 burn-in samples,
respectively. The 5000 sweeps after the burn-in period are thinned with
50 draws to reduce dependence of iterations. The first column of Figure~\ref{fig:normexhier_s} shows
the probabilities of regime indicator $ s_t = i $ of the two
models. Intersections of the lines $ s_t = i $
clearly demonstrate the break locations.

To compare our proposed DPHMM to Chib's method, we also report the posterior
inference of Chib's model under the true change-point number and the same
model specification as in (\ref{eq:norm_theta_post}) and
(\ref{eq:sig2_post}).\footnote{The prior of the transition probabilities in \citet{chib98}'s model
is assumed to be $\mathrm{Beta}(a,b)$. The parameters $a$ and $b$ are chosen
to reflect equidistant duration of each state. For example, in the case of one
change-point
with $n = 150$ sample size, we take $b = 0.1$ and $a = n/2 \times b = 7.5$,
i.e. $\mathrm{Beta}(7.5,0.1)$.}
The posterior means and standard deviations of parameters are summarized in
Table~\ref{tab:post_normodel1_2}.
First, our method performs well in all cases where the posterior distributions
concentrate on the true values. The sample first-order serial
correlations demonstrate good mixing of the samplers.
Second, our results are comparable to those estimated from Chib's model.
The Bayes factors show that in most of the cases the models with the
true number of change-points are preferred to others. For example, in Model 2
where two change-points exist, the Bayes factors comparing models with $k=1$ versus $k=2$ are
close to zero favoring the two change-point model. Likewise, the Bayes
factors comparing $k=2$ versus $k=3$ favor the two change-point model.
Hence we conclude that the model with two
change-points is correctly specified with high probability.
However, the Bayes factor fails to detect the correct number of changes in
Model 1 with unknown variance. The values suggest a model with two change
points. The posterior probabilities of states
estimated by Chib's model are shown in the second column of
Figure~\ref{fig:normexhier_s}.
In summary, the simulation results demonstrate that our algorithm
works well in the normal mean-shift models and is robust to the change-point
number compared to Chib's model.

\begin{table}[t!]
\caption{Posterior estimates of Normal Mean-Shift Models 1 and 2. Mean and SD
  denote, respectively, posterior mean and posterior
  standard deviation.}
\label{tab:post_normodel1_2}
\begin{tabular*}{\columnwidth}{@{\extracolsep{\fill}}cccccccc}
  \hline
  & \multicolumn{7}{c}{\bf Model 1} \\
  & \multicolumn{3}{c}{Known variance} &&  \multicolumn{3}{c}{Unknown variance}\\
  & Mean & SD & True value && Mean & SD & True value \\
  \hline
  & \multicolumn{7}{c}{\textit{DPHMM}}
  \\
  $ \theta_1 $ & 0.9123 & 0.2519 & 1.0000   && 0.9154 & 0.2432 & 1.0000\\
  $ \theta_2 $ & 2.9375 & 0.1812 & 3.0000   && 2.9354 & 0.1723 & 3.0000\\
  $ \sigma^2 $ &        &        &          && 2.8244 & 0.3343 & 3.0000\\
  \\
  & \multicolumn{7}{c}{\textit{Chib's Model with $k=1$}}
  \\
  $ \theta_1 $ & 0.8980 & 0.1889 & 1.0000 && 0.9521 & 0.2441 & 1.0000 \\
  $ \theta_2 $ & 2.9520 & 0.1324 & 3.0000 && 2.9331 & 0.1619 & 3.0000 \\
  $ \sigma^2 $ &        &        &   && 2.7933 & 0.4121 & 3.0000 \\
  & \multicolumn{7}{c}{\textit{Bayes Factor Analysis}} \\
  $ k = 1 $ vs. $ k = 2$ & \multicolumn{2}{c}{1.730} & && \multicolumn{2}{c}{0.548} &\\
  $ k = 1 $ vs. $ k = 3$ & \multicolumn{2}{c}{1.374} & && \multicolumn{2}{c}{1.010} &\\
  $ k = 2 $ vs. $ k = 3$ & \multicolumn{2}{c}{0.794} & && \multicolumn{2}{c}{1.850} &\\
  \hline
  & \multicolumn{7}{c}{\bf Model 2} \\
  & \multicolumn{3}{c}{Known variance} &&  \multicolumn{3}{c}{Unknown variance}\\
  & Mean & SD & True value && Mean & SD & True value \\
  \hline
  & \multicolumn{7}{c}{\textit{DPHMM}}
  \\
  $ \theta_1 $ & 1.1746 & 0.2777 & 1.0000 && 1.2770 & 0.2584 & 1.0000\\
  $ \theta_2 $ & 2.9758 & 0.2874 & 3.0000 && 3.2176 & 0.2672 & 3.0000\\
  $ \theta_3 $ & 5.3344 & 0.2459 & 5.0000 && 5.1108 & 0.2682 & 5.0000\\
  $ \sigma^2 $ &        &        &        && 3.1102 & 0.3726 & 3.0000\\
  \\
  & \multicolumn{7}{c}{\textit{Chib's Model with $k=2$}}
  \\
  $ \theta_1 $ & 1.1770 & 0.2052 & 1.0000 && 1.3680 & 0.2643 & 1.0000\\
  $ \theta_2 $ & 2.9480 & 0.1916 & 3.0000 && 3.1920 & 0.2730 & 3.0000\\
  $ \theta_3 $ & 5.3320 & 0.1823 & 5.0000 && 5.0310 & 0.2461 & 5.0000\\
  $ \sigma^2 $ &        &        &   && 3.0780 & 0.6794 & 3.0000\\
  & \multicolumn{7}{c}{\textit{Bayes Factor Analysis}} \\
  $ k = 1 $ vs. $ k = 2$ & \multicolumn{2}{c}{0.000} & && \multicolumn{2}{c}{0.0193} &\\
  $ k = 1 $ vs. $ k = 3$ & \multicolumn{2}{c}{0.000} & && \multicolumn{2}{c}{0.0583} &\\
  $ k = 2 $ vs. $ k = 3$ & \multicolumn{2}{c}{3.090} & && \multicolumn{2}{c}{3.0332} &\\
  \hline
\end{tabular*}
\end{table}

\subsection{Robustness Check of Change-Point Number}
In this section, we study the robustness of our algorithm in detecting the
true number of change-points. Although our
method does not require prespecification of the change-point number, it is
still possible that our algorithm fails to estimate the correct number of
change-points. Thus,
we replicate the entire estimation process as in the previous section for\ 1000
times and record the estimated change-point number in each replication. Specifically,
in each of the 1000 replications, we iterate the Gibbs sampler for 5000 times
and the change-point number of the last sample is recorded.
Therefore, we obtain 1000 collections of change-point numbers.
Table \ref{tab:failing_rate} reports the frequencies of the detected change-point numbers.
We see that with high frequency (over 99\%) our method detects one change-point ($ k = 1 $)
in Model 1 in both cases of known and unknown variance. In Model 2,
our results show that over 90\% of the 1000 replications detect two
change-points ($k = 2$), and over 99\% detect at least one
change-point. The\ figures demonstrate that our
algorithm correctly detects the change-point
number with high probability in different cases.

\begin{table}[htbp]\centering
\caption{Frequencies of estimated change-point numbers.}
\label{tab:failing_rate}
\begin{tabular}{lrrrrrrr}
  \hline
  & \multicolumn{3}{c}{Known variance} && \multicolumn{3}{c}{Unknown variance} \\
  & $k = 0$ & $k = 1$  & $k = 2$ && $k = 0$ & $k = 1$ & $k = 2$ \\
  \hline
      {\bf Model 1}      & \multirow{2}{*}{0.3\%} & \multirow{2}{*}{99.7\%}  & \multirow{2}{*}{0.0\%} &&
      \multirow{2}{*}{0.5\%} & \multirow{2}{*}{99.5\%} & \multirow{2}{*}{0.0\%} \\
      (One change-point) & &&   &&   &    \\
      {\bf Model 2}      & \multirow{2}{*}{0.0\%} & \multirow{2}{*}{6.5\%}  & \multirow{2}{*}{93.5\%} &&
      \multirow{2}{*}{0.2\%} & \multirow{2}{*}{8.7\%} & \multirow{2}{*}{91.1\%} \\
      (Two change-points)& &&   &&   &       \\
      \hline
\end{tabular}
\end{table}

\section{Learning $\alpha$ and $\beta$}\label{sec:learn_albe}
In the previous section, we set the DPHMM parameters $\alpha = 3$ and $\beta =
2$ in estimating the simulated models. In order to learn about $\alpha $ and $\beta
$, we propose to use vague Gamma priors, see \citet{beal2002}.
Note that with the number of states specified in each MCMC sweep,
the DPHMM reduces to the generalized Dirichlet distribution (GDD),
see \citet{connor_mosimann1969} and \citet{wong_GDD_1998}.
Hence the posterior is
\begin{equation}\label{eq:alpha_beta_post}
p(\alpha,\beta\mid S_n) \propto \mathrm{Gamma}(a_\alpha,b_\alpha)
                   \mathrm{Gamma}(a_\beta,b_\beta)
                   \prod_{i=1}^{k+1}\frac{\beta \Gamma(\alpha+\beta)}{\Gamma(\alpha)}
                   \frac{\Gamma({n_{ii}+\alpha})}{\Gamma({n_{ii}+1+\alpha+\beta)}},
\end{equation}
where $ n_{ii} = \sum_{t=1}^{T-1}\delta(s_{t}, i)\delta(s_{t+1},i) $
denotes the counts of self transitions. We set $a_\alpha=b_\alpha=a_\beta =
b_\beta = 1$ here and in the subsequent sections. Below we consider two
alternative approaches for sampling: the first based on maximum-a-posteriori
(MAP) estimation and a second approach using a random walk sampler.

\subsection{The Maximum-a-Posteriori}
We first solve for the maximum-a-posteriori (MAP) estimates
for $\alpha $ and $\beta$ which are obtained as the solutions to the following
gradients using the Newton-Raphson method,
\begin{equation}\label{eq:gradient}
  \begin{aligned}
  \frac{\partial \ln p(\alpha,\beta\mid S_n)}{\partial \alpha}
  &= \frac{a_\alpha-1}{\alpha} - b_\alpha \\
  &+ \sum_{i=1}^{k+1} \left[ \psi(\alpha+\beta) + \psi(n_{ii} + \alpha)
  - \psi(\alpha) - \psi(n_{ii}+1+\alpha+\beta) \right] = 0\\
  \frac{\partial \ln p(\alpha,\beta\mid S_n)}{\partial \beta}
  &= \frac{a_\beta-1}{\beta}-b_\beta + \sum_{i=1}^{k+1}\left[
  \frac{1}{\beta} + \psi(\alpha+\beta)
  - \psi(n_{ii}+1+\alpha+\beta)  \right] = 0,
  \end{aligned}
\end{equation}
where $\psi(\cdot)$ is the digamma function defined as
$\psi(x) = d \ln \Gamma(x)/ dx$.

We implement our algorithm in the previous section together with
the MAP update of $\alpha$ and $\beta$ in each sweep.
The DPHMM with MAP update correctly detects the true number of change-points
in all cases.
Table \ref{tab:MAP} shows the MAP solutions for $\alpha$ and $\beta$, and
the posterior estimates of all parameters in each model. We can see
that the average MAP values of $\alpha$ and $\beta$ are $ 0.6353$ and $0.1937$
respectively in Model 1 with known and unknown variance. The results
are slightly different in Model 2 such that average MAP values
are $ 0.9451$ and $0.2360$ respectively. We also report the sample standard
errors which show evidence of stability of the MAP values after the
burn-in period. In all cases, $\alpha$ is greater than $\beta$ indicating that
the algorithm tends to linger in existing states rather than exploring a new
one. Besides, all parameter estimates are in line with the results in the
previous section when $\alpha$ and $\beta$ are prespecified.

\begin{table}[t]
  \caption{MAP of $\alpha$ and $\beta$ in Normal Mean-Shift Models 1 and
      2. Average MAP values of $\alpha$ and $\beta$ are reported with
      standard deviations within parentheses. For other parameters, the values
      are posterior means and posterior standard deviations.}
  \label{tab:MAP}
    \begin{tabular*}{\columnwidth}{@{\extracolsep{\fill}}lcccc}
      \hline
      & \multicolumn{2}{c}{\bf Model 1} & \multicolumn{2}{c}{\bf Model 2}\\
      & \textit{known variance} &  \textit{unknown variance}
      & \textit{known variance} &  \textit{unknown variance}\\
      \hline
      $ \alpha $   &  0.6353 (0.0007) &  0.6353 (0.0007)  & 0.9453 (0.0008) & 0.9451 (0.0005)\\
      $ \beta $    &  0.1937 (0.0005) &  0.1937 (0.0005)  & 0.2364 (0.0001) & 0.2361 (0.0005)\\
\\
      $ \theta_1 $ &  0.9145 (0.2560) &  0.9049 (0.2430)  & 1.1843 (0.2847) & 1.2675 (0.2585)\\
      $ \theta_2 $ &  2.9326 (0.1745) &  2.9385 (0.1731)  & 2.9859 (0.2892) & 3.2222 (0.2716)\\
      $ \theta_3 $ &                  &                   & 5.3349 (0.2468) & 5.1026 (0.2661)\\
      $ \sigma^2 $ &                  &  2.8207 (0.3340)  &                 & 3.0908 (0.3752)\\
      \hline
    \end{tabular*}
\end{table}

\subsection{The Metropolis-Hastings Sampler}\label{sec:MH}
We also consider a Metropolis-Hastings (M-H) sampler for the
posterior (\ref{eq:alpha_beta_post}). The candidate-generating density
is assumed to be the
random walk process with positive support
\begin{equation}
  f(\alpha'|\alpha) \propto \phi(\alpha'-\alpha),\quad  \alpha'> 0,
\end{equation}
where $\alpha$ is the value of the previous draw, $\phi(\cdot)$ is the
standard normal density function. The acceptance ratio given $\beta$ is thus
\begin{equation}
  A(\alpha,\alpha') = \frac{p(\alpha',\beta\mid S_n) \;
  \Phi(\alpha)}{p(\alpha,\beta\mid S_n)\;\Phi(\alpha')},
\end{equation}
where $\Phi(\cdot)$ is the standard normal distribution function. The same M-H
sampler is also applied to $\beta$ given the updated
$\alpha$. We incorporate the M-H sampler of $\alpha$ and $\beta$ in the Gibbs
sampler in Section 5. The posterior estimators are
shown in Table \ref{tab:MH}.
The posterior means and standard deviations of parameter $\theta_i$ and $\sigma^2$ are
similar to those obtained in the previous analyses. The posterior mean of $\alpha$ is greater
than the posterior mean of $\beta$ in all models. This affirms the conclusion in the
MAP results that the algorithm tends to linger in existing states
rather than exploring a new one. Both the MAP and M-H methods correctly estimate
the number of change-points in each simulation.

\subsection{Comparison between MAP and M-H}
We can see that the estimates of $\alpha$ and $\beta$ from the two approaches are quite
different as shown in Tables \ref{tab:MAP} and \ref{tab:MH}.
The MAP as a point estimator may not reflect the variations of $\alpha$ and
$\beta$, whereas the M-H is a typical Bayesian method which can be
incorporated into the MCMC sampler of other parameters in question.
Moreover, the MAP approach may be limited when the posterior happens to be multi-modal.
Therefore, the M-H method is preferred in practice and the following empirical
studies are conducted with the M-H sampler.

\begin{table}[t]
  \caption{M-H sampler of $\alpha$ and $\beta$ in Normal Mean-Shift Models 1
  and 2. Posterior means and posterior standard deviations within parentheses.}
  \label{tab:MH}
    \begin{tabular*}{\columnwidth}{@{\extracolsep{\fill}}lcccc}
      \hline
      & \multicolumn{2}{c}{\bf Model 1} & \multicolumn{2}{c}{\bf Model 2}\\
      & \textit{known variance} &  \textit{unknown variance}
      & \textit{known variance} &  \textit{unknown variance}\\
      \hline
      $ \alpha $   &  1.8073 (1.3345) & 1.8145 (1.3646)   &  2.2007 (1.5179) & 2.1482 (1.4731)\\
      $ \beta $    &  0.3446 (0.2168) & 0.3455 (0.2176)   &  0.3667 (0.2087) & 0.3622 (0.2086)\\
      \\
      $ \theta_1 $ &  1.1182 (0.2521) & 1.0839 (0.2399)   &  1.1256 (0.2558) & 1.0892 (0.2457)\\
      $ \theta_2 $ &  2.9447 (0.1724) & 3.0939 (0.1701)   &  2.8277 (0.2644) & 3.0341 (0.2449)\\
      $ \theta_3 $ &                  &                   &  5.0138 (0.2461) & 5.1172 (0.2454)\\
      $ \sigma^2 $ &                  & 2.8190 (0.3298)   &                  & 2.8793 (0.3455)\\
      \hline
    \end{tabular*}
\end{table}

\section{Empirical Applications}\label{sec:app}

\subsection{Poisson Data with Change-Point}

We first apply our Dirichlet process multiple change-point model to
the much analyzed data set on the number of
coal-mining disasters by year in Britain over
the period 1951--1962 (\citet{jarrett1979}, \citet{carlin_gelfand_smith92}
and \citet{chib98}).

\begin{figure}[t]
\includegraphics{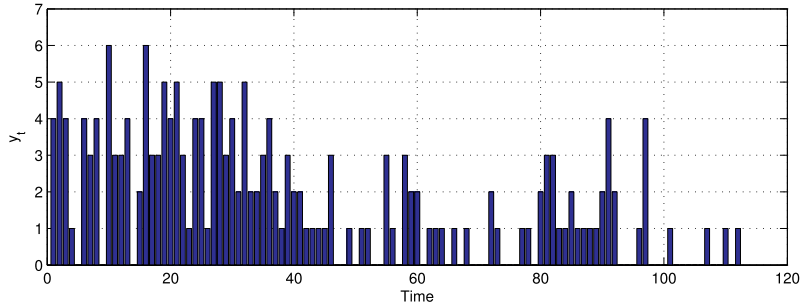}
\caption{Data on coal mining disaster count $ y_t $.}
\label{fig:jarret79coal}
\end{figure}

Let the disaster count $ y $ be modeled by a Poisson distribution
\begin{equation}
  f(y\mid \lambda) = \lambda^y e^{-\lambda}/y!.
\end{equation}
The observation sequence $ Y_n = ( y_1, y_2, \dots, y_{112} )' $
is subject to some unknown change-points. We plot the data $ y_t $ in
Figure~\ref{fig:jarret79coal}.
\citet{chib98} estimates the models with one change-point ($ k = 1$) and with
two change-points ($ k = 2 $), respectively.
He assumes the parameter $ \lambda $ following the prior $ \mathrm{Gamma}(2,1) $
in the one-change-point case and the prior $ \mathrm{Gamma}(3,1) $ in the other.
Hence, given the regime indicators $ S_n $,
the corresponding parameter $ \lambda_i $ in regime $ i $ has the following posteriors
with respect to the two priors:
\begin{equation}\label{eq:post1_lambda}
  \mbox{Posterior 1:}\;\;\;\;
\lambda_i\mid S_n, Y_n \sim \mathrm{Gamma}(\lambda_i\mid 2 + U_i, 1 + N_i),
\end{equation}
and
\begin{equation}\label{eq:post2_lambda}
  \mbox{Posterior 2:}\;\;\;\;
\lambda_i\mid  S_n, Y_n \sim \mathrm{Gamma}(\lambda_i\mid 3 + U_i, 1 + N_i),
\end{equation}
where $ U_i = \sum_{t=1}^{112} \delta(s_t, i)y_t $ and $ N_i =
\sum_{t=1}^{112} \delta(s_t, i) $. We perform our algorithm  with the
following Gibbs steps:

\vspace{3 mm}
\emph{Step} 1. Sample $ S_n \mid \lambda, Y_n $
as in (\ref{eq:gibbs_post_state}) and obtain $ k $,

\vspace{3 mm}
\emph{Step} 2. Sample $ \lambda_i \mid S_n, Y_n $
as in (\ref{eq:post1_lambda}) or (\ref{eq:post2_lambda}),
\vspace{3 mm}

\emph{Step} 3. Update $ \alpha $ and $\beta$ with the Metropolis-Hastings Sampler
as in Section \ref{sec:MH}.
\vspace{3 mm}

\begin{figure}[t!]
\includegraphics{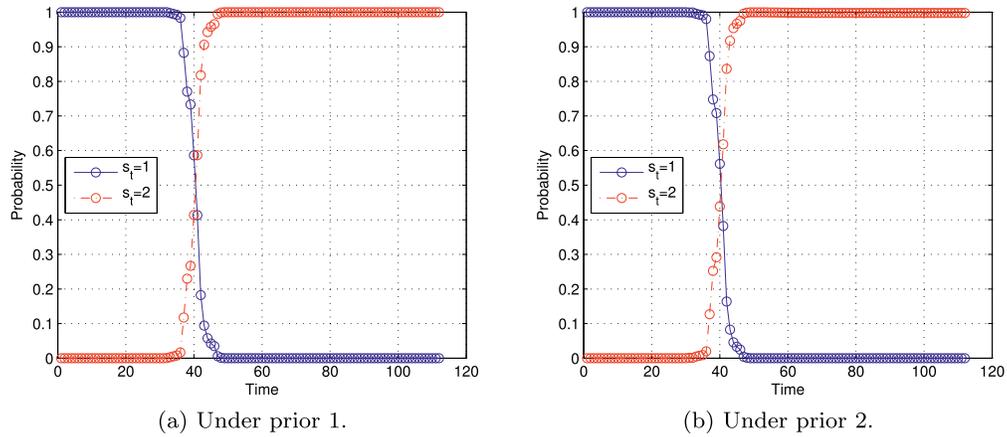}
\caption{Posterior probability of $ s_t = i $.}
\label{fig:poisson_s}
\end{figure}

\begin{figure}[t]
\centering
\includegraphics{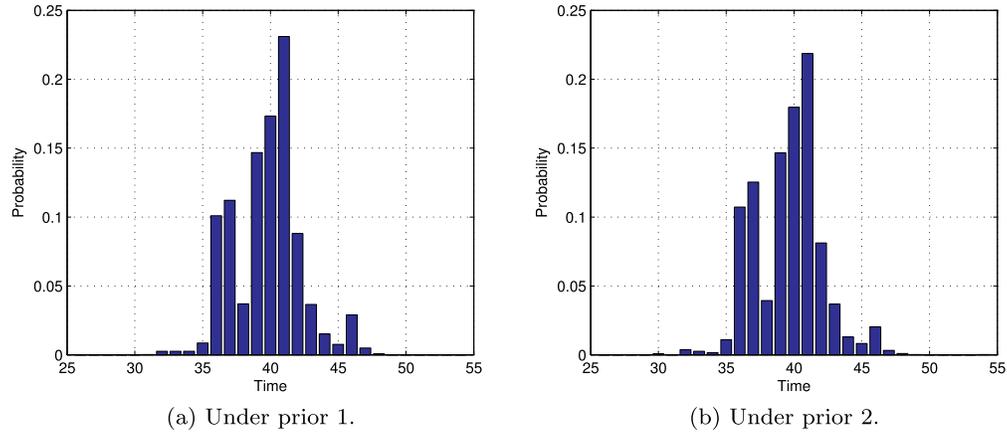}
\caption{Posterior probability mass function of change-point location $ \tau_i $.}
\label{fig:poisson_post}
\end{figure}

The above Gibbs sampler is conducted for 5000 sweeps
with 1000 burn-in samples. To reduce the sampler dependency, the 5000
sweeps are thinned by 50 draws. The sampler estimates one change-point in the data.
Figure~\ref{fig:poisson_s} shows the posterior
probabilities of the regime indicator $ s_t = i $ at each time point $ t $. The
intersections of the two lines $ s_t = 1 $ and $ s_t = 2 $ show that
the break location exits at around $ t = 40 $. Figure~\ref{fig:poisson_post}
provides the distribution of the transition points $ \tau_i $. Interestingly,
our model produces exactly the same figure as the one in \citet{chib98}. The
change-point is identified as occurring at around $ t = 41 $.

The corresponding posterior
means of the parameters $ \lambda_1 $ and $ \lambda_2 $ are
$ 3.1006 $ and $ 0.9387 $ with posterior
standard deviations, $ 0.2833 $ and $ 0.1168 $, respectively, under the
prior $ \mathrm{Gamma}(2,1) $. The posterior means of $\alpha$ and $\beta$ are $1.8101$ and $0.3697$ with
standard deviations $1.3577$ and $0.2464$.
When using the prior $ \mathrm{Gamma}(3,1) $,
we have the posterior means of $ \lambda_1 $ and $ \lambda_2 $ equal to
$ 3.1308 $ and $ 0.9567 $ with posterior
standard deviations~$ 0.2877 $ and $ 0.1218 $, respectively.
The posterior means of $\alpha$ and $\beta$ are $1.8375$ and $0.3715$ with
standard deviations $1.3456$ and $0.2360$, respectively under
prior 2.
All our results closely match those of the literature and we show a
certain robustness of our model under different prior assumptions.

In order to check the robustness of the estimation of the number of change-points, $ k $,
we conduct 1000 replications of the above estimation process and collect 1000
change-point numbers.
When the first prior is assumed, $ 77.23\% $ of the 1000 replications detect one
change-point. We find a similar
result for prior 2. Hence, we conclude that without assuming the number
of change-points \textit{a priori}, our algorithm detects the same
change-point number as in the model developed by \citet{chib98}
with high probability.

\subsection{Real Output}

\begin{figure}[t]
\includegraphics{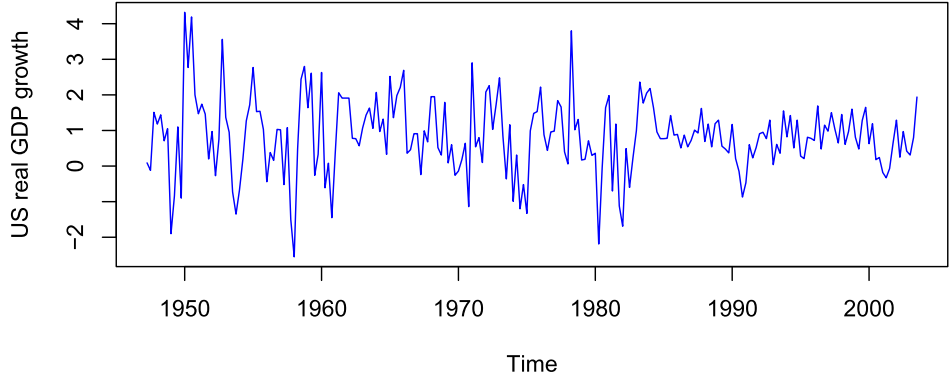}
\caption{US real GDP growth from the second quarter of 1947 to the third
quarter of 2003.}
\label{fig:usgdp}
\end{figure}

We also apply our algorithm to estimate structural changes in real Gross
Domestic Product growth. The data and model are drawn from \citet{maheu_gordon08} (see also
\citet{Geweke_yu2011}). Let $ y_t =
100[\log(q_t/q_{t-1})-\log(p_t/p_{t-1})] $, where $ q_t $ is quarterly US GDP
seasonally adjusted and $ p_t $ is the GDP price index. The data range from
the second quarter of 1947 to the third quarter of 2003, for a total of 226 observations (see Figure~\ref{fig:usgdp}).
We model the data with a Bayesian \textsc{ar}$(2)$ model with
structural change. The frequentist autoregressive structural change-model can be found
in \citet{chong2001}. Suppose the data are subject to $ k $ change-points and
follow
\begin{equation}
  y_t = \beta_{0,s_t} + \beta_{1,s_t} y_{t-1} + \beta_{2,s_t} y_{t-2} +
  \varepsilon_t, \;\;\;\;\varepsilon_t \sim \mathrm{N}(0,
  \sigma^2_{s_t}),\;\;\;\;s_t = 1, 2, \dots, k+1.
\end{equation}
We assume the following hierarchical priors to $ \beta_{0,i},
\beta_{1,i} $ and $\beta_{2,i} $:
\begin{equation}
  \beta_i = (\beta_{0,i}, \beta_{1,i}, \beta_{2,i})' \sim
  \mathrm{N}(\mu, V), \;\;\;\; i = 1, \dots,
  k+1,
\end{equation}
where $ \mu = ( \mu_0, \mu_1, \mu_2 )' $ and $ V =
\mathrm{Diag} ( v^2_0, v^2_1, v^2_2 ) $, such that
\begin{equation}
  p(\mu_j, v_j^2) \propto \mbox{Inv-Gamma}(v_j^2\mid a, b),\;\;\;\;j = 0,1,2.
\end{equation}
We assume the noninformative prior for $ \sigma^2_i $ such that
\begin{equation}
    p(\sigma_i^2) \propto 1/\sigma_i^2, \;\;\;\;
    i = 1, \dots, k+1.
\end{equation}

Conditional on $ \sigma_i^2 $, the sampling of $ \beta_i $, $
\mu $ and
$ V $ is similar to Section \ref{sec:normex}. For
$ \sigma_i^2 $, we can draw from the following full conditional:
\begin{equation}
  \sigma_i^2\mid \beta_i,S_n, Y_n \sim
  \mbox{Inv-}\chi^2\left(\sigma_i^2 \mid  \tau_i - \tau_{i-1},
  \frac{\omega_i^2}{\tau_i - \tau_{i-1}}\right),
  \;\;\;\; i = 1, \dots, k+1,
\end{equation}
where $ \omega_i^2 = \sum_{\tau_{i-1} < t \leq \tau_i}(y_t - \beta_{0,s_t} -
\beta_{1,s_t} y_{t-1} - \beta_{2,s_t} y_{t-2})^2$.

\begin{figure}[t]
\includegraphics{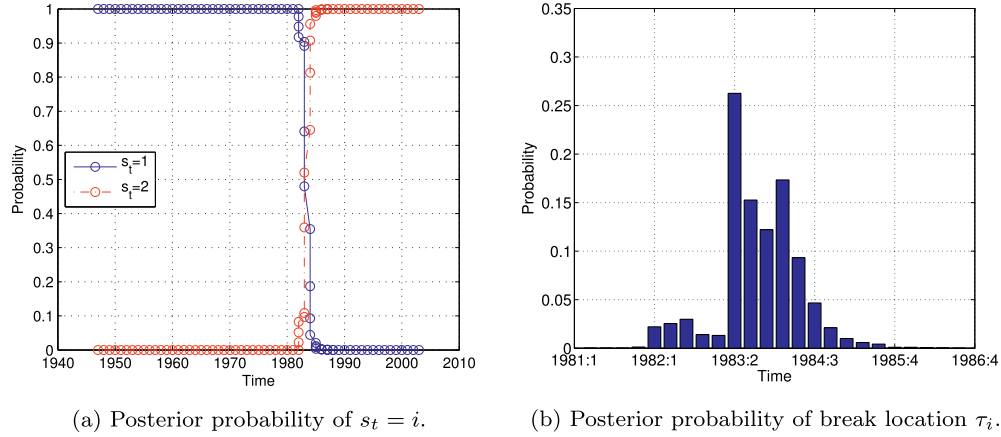}
\caption{US real GDP growth structural change model.}
\label{fig:usgdp_ks}
\end{figure}

As in the previous applications, we set the inverse-Gamma hyperparameters\break \mbox{$ a = b =1 $}.
The M-H update of $\alpha$ and $\beta$ follows the discussion in
Section \ref{sec:MH}. The Gibbs sampler is conducted for 5000 sweeps with 1000
burn-in samples. The 5000 sweeps are thinned by 50 draws.
The posterior
probabilities of the regime indicator $ s_t $ in Figure~\ref{fig:usgdp_ks}a suggest that the structural
break exists between the years 1980 and 1990. Figure~\ref{fig:usgdp_ks}b
further shows the change-point at the second quarter of 1983, which is close to the results in
\citet{maheu_gordon08}.\footnote{See Figure 4 in \citet{maheu_gordon08}.} The
posterior estimates are summarized in
Table~\ref{tab:usgdp_postinfer}. Finally, the posterior means of $\alpha $
and $\beta$ are $1.7749$ and $0.3045$ with standard deviations
$ 1.3422 $ and $0.1939$ respectively. All of our results are consistent with Chib's
estimates.

\begin{table}[t]
\caption{US real GDP growth structural change model with one
  change-point. Posterior means and posterior standard deviations within
  parentheses. The results applying Chib's model are drawn
  from \citet{maheu_gordon08}.}
\label{tab:usgdp_postinfer}
\begin{tabular*}{\columnwidth}{@{\extracolsep{\fill}}lcccccc}
  \hline
  && \multicolumn{2}{c}{\bf Chib's Model} && \multicolumn{2}{c}{\bf DPHMM}\\
  && $ s_t = 1 $ & $ s_t = 2 $
  && $ s_t = 1 $ & $ s_t = 2 $ \\
  \hline
  $ \beta_{0, s_t} $ && 0.5642 (0.1228) & 0.4434 (0.1162) && 0.5499 (0.1303) & 0.3894 (0.1169) \\
  $ \beta_{1, s_t} $ && 0.2716 (0.0734) & 0.2792 (0.1052) && 0.2812 (0.0837) & 0.2796 (0.1173)\\
  $ \beta_{2, s_t} $ && 0.0800 (0.0739) & 0.1588 (0.1010) && 0.0913 (0.0855) & 0.2253 (0.1124)\\
  $ \sigma^2_{s_t} $ && 1.3331 (0.1542) & 0.3362 (0.0516) && 1.4089 (0.1722) & 0.2672 (0.0460)\\
  \hline
\end{tabular*}
\end{table}

Finally, we replicate 1000 times the whole estimation process and check the
robustness of the detected change-point number. The result suggests
that nearly 100\% of the replications detect one break point.

\section{Concluding Remarks}\label{sec:rm}

In this paper, we have proposed a new Bayesian multiple change-point
model, that is the Dirichlet process hidden Markov model.
Our model is semiparametric in the sense that the number of states is not built-in to the
model but endogenously determined. As a result, our model avoids the model
misspecification problem.
We have proposed an MCMC sampler which only needs to sample the states around change-points.
We have also proposed the MAP and M-H updates of hyperparameters
in the DPHMM process. We have presented three specific models, namely, the
discrete Poisson model, the continuous normal model, and the \textsc{ar}$(2)$
model with structural changes. Results from the simulations and empirical
applications showed that our Dirichlet process hidden Markov multiple change-point model
detected the true change-point numbers and locations with high
accuracy.\vspace*{-6pt}

\appendix
\renewcommand{\theequation}{A-\arabic{equation}}
\setcounter{equation}{0}

\section*{Appendix}\vspace*{-3pt}

In the appendix, we give the derivations of the full conditionals and
the Gibbs samplers in Section \ref{sec:normex}. For the case of known variance,
we first rewrite the hierarchical model
(\ref{eq:norm_model}) as the joint distribution
\begin{equation}\label{eq:norm_joint}
  p(Y_n, \theta, \mu, \upsilon^2\mid  S_n, \sigma^2) \propto
  \prod_{i=1}^{k+1} \mathrm{N}(\tilde{y}_i\mid  \theta_i, \sigma_i^2)
  \prod_{i=1}^{k+1} \mathrm{N}(\theta_i\mid  \mu, \upsilon^2)
  p(\mu,\upsilon^2),
\end{equation}
where $p(\mu,\upsilon^2)$ corresponds to $\mbox{Inv-Gamma}(\upsilon^2\mid a,
b)$, $ \theta = ( \theta_1, \dots, \theta_{k+1} )' $ and
\begin{equation}\label{eq:repara}
  \tilde{y}_i = \frac{\sum_{\tau_{i-1} < t \leq \tau_{i}} y_t}
        {\tau_{i}- \tau_{i-1}}, \;\;\;\; \mbox{and}\;\;\;\;
        \sigma_i^2 = \frac{\sigma^2}{\tau_{i}-\tau_{i-1}}.
\end{equation}
From (\ref{eq:norm_joint}) and (\ref{eq:repara}), we have the following full conditionals:
\begin{equation}\label{eq:norm_theta_post}  \begin{aligned}
    p(\theta_i\mid \mu, \upsilon^2, \sigma^2, S_n, Y_n) &\propto \prod_{i=1}^{k+1}
    \mathrm{N}(\tilde{y}_i\mid  \theta_i, \sigma_i^2) \mathrm{N}(\theta_i\mid  \mu,
    \upsilon^2)\\
    &\propto \mathrm{N}\left( \theta_i\mid \frac{\tilde{y}_i/\sigma_i^2 +
      \mu/\upsilon^2}{1/\sigma_i^2 + 1/\upsilon^2}, \frac{1}{1/\sigma_i^2 +
      1/\upsilon^2}\right),\\
    p(\mu\mid \theta, \upsilon^2, S_n, Y_n) &\propto \prod_{i=1}^{k+1}
    \mathrm{N}(\theta_i\mid \mu, \upsilon^2) p(\mu, \upsilon^2)\\
    &\propto \mathrm{N}(\mu\mid \bar{\theta}, \upsilon^2/(k+1)),\\
    p(\upsilon^2\mid \theta, \mu, S_n, Y_n) &\propto
    (\upsilon^2)^{-(k+1)/2}\exp\left\{ -\frac{1}{2}\sum_{i=1}^{k+1}(\theta_i -
    \mu)^2/\upsilon^2\right\} p(\mu, \upsilon^2)\\
    &\propto \mbox{Inv-Gamma}\left(\upsilon^2\mid a + \frac{k+1}{2},
    b + \frac{1}{2}\sum_{i=1}^{k+1}(\theta_i-\mu)^2\right),\\
  \end{aligned}
\end{equation}
where $ \bar{\theta} = \sum_{i=1}^{k+1}\theta_i / (k+1)$. Therefore,
we can perform the following Gibbs sampler:

\vspace{2 mm}
\indent \emph{Step} 1. Sample $ S_n \mid \theta, \mu, \upsilon,
  Y_n $ as in (\ref{eq:gibbs_post_state}) and obtain $ k $,

\indent \emph{Step} 2. Sample $ \theta, \mu, \upsilon\mid S_n,
  Y_n $ as in (\ref{eq:norm_theta_post}).
\vspace{2 mm}

\noindent For the case of unknown variance, the full conditional with respect to
(\ref{eq:sig2_prior}) is
\begin{equation}\label{eq:sig2_post}
  \sigma^2\mid \theta,S_n, Y_n \sim
  \mbox{Inv-Gamma}\left(\sigma^2 \bigg\vert  c + \frac{n}{2},
    d + \frac{1}{2}\sum_{t=1}^n(y_t-\theta_{s_t})^2\right).
\end{equation}

\noindent Conditional on $ \sigma^2$, we apply the same estimation strategy
discussed above. The Gibbs sampler is thus

\vspace{2 mm}
\emph{Step} 1. Sample $ S_n \mid  \theta, \mu, \upsilon,
  \sigma^2, Y_n $ as in (\ref{eq:gibbs_post_state}) and obtain $ k $,

\emph{Step} 2. Sample $ \theta, \mu, \upsilon, \sigma^2 \mid S_n,
  Y_n $ as in (\ref{eq:norm_theta_post}) and (\ref{eq:sig2_post}).

\begin{acknowledgement}
The authors would like to thank all the participants in the Econometric
Society Australasian Meeting 2011, Adelaide, Australia, July 2011, for
helpful comments and discussions.

The third author (P.G.) acknowledges the support of DST grant
(SR/S4/MS:648/10) from the Government of India.

\end{acknowledgement}

\end{document}